\documentclass[preprint,12pt]{elsarticle}
\def\R{\mathbb{R}}
\def\Z{\mathbb{Z}}

\def\({\left(}
\def\){\right)}

\def\Re{\hbox{Re }}
\def\Im{\hbox{Im \  }}

\def\ep{\varepsilon}

\usepackage{amsmath, graphicx,subfigure,cases,amssymb,amsfonts,amsthm}

\newtheorem{theo}{Theorem}

\newtheorem{prop}[theo]{Proposition}
\newtheorem{rem}{Remark}

\usepackage{hyperref}

\begin{document}
\begin{frontmatter}

\title{Solitons and solitonic vortices in a strip}

\author{Amandine Aftalion }

\address{Ecole des Hautes Etudes en Sciences Sociales,  Centre d'Analyse et de Math\'ematique Sociales, UMR-8557, Paris, France.}
\ead{amandine.aftalion@ehess.fr}

\author{Etienne Sandier}
\address{LAMA, Univ Paris Est Creteil, Univ Gustave Eiffel, UPEM, CNRS, F-94010 Créteil, France.}
\ead{sandier@u-pec.fr}
\vspace{10pt}

\begin{abstract}
We study the ground state of the Gross Pitaevskii energy in a strip, with a phase imprinting  condition, motivated by recent experiments on matter waves solitons. We prove that when the width of the strip is small, the ground state is a one dimensional soliton. On the other hand, when the width  is large, the ground state is a solitonic vortex. We provide an explicit expression for the limiting phase of the solitonic vortex as the size of the strip is large: it has the same behaviour as the soliton in the infinite direction and decays exponentially due to the geometry of the strip, instead of algebraically as vortices in the whole space.
\end{abstract}

%
%
%

%

\end{frontmatter}

\section{Introduction}

The instability of solitons is a phenomenon which has been widely studied both from a mathematical and a physical point of view. A soliton is a one dimensional solitary wave propagating at constant velocity. The simplest situation is actually the case of zero velocity. Solitons have been observed in a large variety of nonlinear media \cite{dauxois,newell}. The instability has been investigated in nonlinear optics, crystals, and more recently in ultra cold atomic vapors and polaritons. However a full characterization of the instability of solitons is still missing.

\subsection{Physical motivation}

Ultra cold atoms can be controlled at will and therefore allow to study collective behaviours. They are described by a wave function whose phase plays the role of orientation in magnetization. A characteristic feature is the observation of defects \cite{chevysol}. The simplest defects in cold atoms are dark solitons which are solitary waves that maintain their shape; they correspond to an envelope having a density dip with a $\pi$ phase shift across its density minimum.  Other types of defects are vortices or vortex rings in which the order parameter winds around a hole or a closed loop where the density vanishes. The study of solitons \cite{gravsmets,bgs,tao}, traveling waves \cite{betsaut,betsaut2,maris,ruiz1,ruiz2} and vortex rings \cite{robrings} in the whole space has been at the center of many mathematical works.

Matter waves dark solitons can be created in experiments by means of various methods \cite{frantz}: phase imprinting, density engineering, quantum state engineering which is a combination of the previous two, matter wave interference and dragging an obstacle through a condensate. Several experimental groups attempted to study solitons by imposing a phase shift in an elongated condensate for respectively bosonic atoms, rubidium \cite{PRLsoli} and sodium \cite{kibble,PRLsol2} and  for fermionic atoms (lithium) \cite{heavysol,PRLrings,PRLsol}. From the first picture, they thought they had observed solitons \cite{kibble,heavysol}. Further investigations were needed to fully understand the phenomena: in the case of lithium, they realized it was not a soliton but thought it was a vortex ring \cite{PRLrings}, until \cite{PRLsol} argued that in fact it was a single straight vortex called solitonic vortex; in sodium it was also confirmed it was not a soliton but a solitonic vortex \cite{PRLsol2}.

 A solitonic vortex is a vortex in a channel whose transverse size is much smaller than its length and has the same asymptotic phase profile as a soliton in the transverse direction. Indeed, by virtue of the transverse boundaries, the solitonic vortex is exponentially localized in the longitudinal direction on the length scale of the transverse dimension. This is a remarkable property since the energy density of isolated vortices or of vortex-antivortex pairs decays algebrically in the absence of boundaries \cite{bgs,betsaut,betsaut2,JR}. Moreover, in the absence of boundaries, vortex-antivortex solutions do not exhibit any phase difference at infinity.
 Mathematically, as we will see, the exponential localisation and the phase difference between the two ends of the channel result from the infinite and periodic structure of image vortices. In a finite channel, because the vortex loses its long-range nature beyond the transverse length scale, this has justified the name "solitonic vortex" introduced in \cite{brand,kopapasol,toikka}.

 In the experiments, it may be that the solitonic vortex is due to the decay of a soliton after what is called a snake instability. It was observed that the soliton is stable only at sufficiently low particles number or high aspect ratio \cite{kopapasol}. The mode associated with snake instability is the solitonic vortex \cite{brand}. The metaphor of snake instability has been used to refer to the bending of the solitonic wave front and the decay of the soliton into closed loops of vortices. In a series of studies, Komineas et al \cite{kopapasol,kopapawaves,komi,kopapaL} analyze how the soliton bifurcates according to the channel length. In some intermediate cases, the soliton initially deforms to become a pair of vortex-antivortex or a vortex ring in 3d. This structure eventually decays into a stable solitonic vortex. The numerics reveal that the vortex-antivortex pair or vortex ring is unstable, but it is sufficiently long lived to be observed both in the numerics and the experiments.

 From a mathematical point of view, Rousset and Tzvetkov, \cite{rousset}, Theorem 3.3, have proved the transverse instability of solitons in the whole space. Nevertheless, nothing is known about the mode of destabilization and whether it turns into a solitonic vortex or a pair of vortices in dimension 2, or a vortex ring in dimension 3. Nothing has been analyzed in a strip yet.

 \hfill


\subsection{Main result}

The aim of this paper is to determine the ground state of the Gross Pitaevskii energy when both a phase imprinting  condition and a reduction of density are imposed. In particular, we want to discriminate the cases where the soliton can be a stable minimizer. We will restrict to static solutions, that is for a zero velocity. In the following, we will set the problem in two dimensions.

\hfill

Let $\Omega_d=\R \times (-d,d) $ be a strip. We will consider the width $d$ as a parameter. Let $u(x,y)$ be a minimizer of
\begin{equation}\label{energy}
E_{d}(u)= \int_{\Omega_d}\left (\frac 12 |\nabla u (x,y)|^2+ \frac 1 {4} (1-|u (x,y)|^2)^2 \right )
 \,dx \ dy,
\end{equation}under the conditions
\begin{equation}\bar u (x,y)= u(-x,y) \hbox{ and } u(-x,-y)=-u(x,y).\label{eqsymmetry}\end{equation} Here the bar denotes the complex conjugate.
 The first equation corresponds to what is called the phase imprinting process: one half of the domain is exposed to a blue detuned laser beam that causes a phase shift of the order parameter  \cite{brand,kopapasol}.
 The second condition corresponds to the fact the density is depleted at a desired location \cite{spielman}, here the origin. In order to impose the fact that the wave function vanishes at the origin, we have decided to impose imparity of the wave function.

These two conditions are equivalent to
\begin{equation}\label{eqsym}\Re u (-x,y)= \Re u (x,y),\quad \Im u(-x,y)=-\Im u(x,y),\end{equation}
\begin{equation}\label{eqsym2}\Re u (x,-y)=- \Re u (x,y),\quad \Im u(x,-y)=\Im u(x,y).\end{equation} Because of the finite energy condition, $|u|$ tends to 1 at infinity, and by the parity, imparity conditions, in fact we will show that the behaviour at infinity is
\begin{equation}\lim_{x\to -\infty} u(x,y)=-i \hbox{ and } \lim_{x\to +\infty} u(x,y)=i,\label{boundcond}\end{equation}
for either $u$ or $\bar u$.

Any ground state is a solution of
\begin{equation}-\Delta u + u (1-|u |^2)=0\end{equation} with the  conditions (\ref{eqsymmetry})-(\ref{boundcond}) and
 $\partial u/\partial n=0$ on $y=\pm d$.  There are special solutions to this equation:
 \begin{itemize}
 \item solitons which are independent of $y$ and are thus one dimensional solutions, $i\tanh (x /\sqrt 2)$,
 \item solitonic vortices which have a zero at the origin and a degree 1.
 \end{itemize} We are going to prove that these solutions characterize the ground states.

 \begin{theo}\label{theo1}There exists a minimizer $u$ of (\ref{energy}) under the conditions (\ref{eqsymmetry}). The behaviour at infinity is given by (\ref{boundcond}) either for $u$ or $\bar u$. If $d$ is sufficiently small, the minimizer is independent of $y$ and is the soliton solution $i\tanh (x /\sqrt 2)$.

 If $d$ is sufficiently large, any minimizer  has a vortex of degree $+1$ or $-1$. Taking the complex conjugate if necessary, we may assume the vortex has degree $+1$. If $u_d(x,y)$ denotes such  a minimizer in $\Omega_d=\R \times (-d,d)$, and if we let
 \begin{equation}\label{vd}v_d(x,y) = u_d(dx, dy),\end{equation}
 then $v_d$, defined in $\Omega_1$,  has the following behaviour: \begin{itemize}\item the map $v_d$ converges as $d\to +\infty$, in $C^1(K)$ for any compact subset $K$ of $\bar \Omega_1\setminus\{0\}$, to the map $e^{i\varphi_0}$, where
 \begin{equation} \varphi_0 =\arctan \frac  {\sinh (\frac {\pi x}{2d}) } {\sin (\frac {\pi y}{2d} )}.\end{equation}
 \item as $d\to +\infty$,
 \begin{equation}\label{bbh}
 E_{d}(u_d) = \pi\log d+ \gamma   - \pi \log\frac\pi 4  + o(1). \end{equation}
 where the constant $\gamma$  is defined in \cite{bbh} as
$$  \gamma = \lim_{r\to +\infty}  \(E_1(u_0,B(0,r))  - \pi \log r\),$$
where $u_0$ is  the unique radial solution of  $-\Delta u = u(1-|u|^2) $ in $\R^2$ of degree +1.\end{itemize}
 \end{theo}
 Note that it is equivalent to use the width $d$ as a parameter or to put a coefficient $1/ d^2$ in front of the term $(1-|u|^2)^2$ using the change of function (\ref{vd}).

  The first part of the Theorem is proved using the Maximum Principle in narrow domains which provides uniqueness.
 For the second part, we introduce the stream function $w_0$, which is such that $\varphi_0$ is its imaginary part, and the exponential of the harmonic conjugate of $\varphi_0$ is its modulus. More precisely,
 $h_0 = \log|w_0|$, where  \begin{equation}\label{w00}w_0(z) = i \tanh\(\frac\pi{4} z\)\end{equation} and $h_0$ is the harmonic conjugate of $\varphi_0$. Once the functions for the strip are computed, the proof follows the main ideas of the analysis of vortices in the seminal book by Bethuel, Brezis, Helein \cite{bbh} and the analysis of vortex balls by Sandier \cite{sandier} and Jerrard  \cite{jerrard}.

\subsection{Further questions}

Let us point out that without the second condition (\ref{eqsymmetry}), there is a minimizing sequence whose energy tends to 0, namely, for large $h$, the functions $u_h (x)=e^{\frac{i\pi x}{2 h}}$ which are defined in $(-h,h)$ and extended by $i$ and $-i$ outside.  Therefore, to have the existence of a minimizer, we have to impose extra conditions, and in particular one corresponding to a density depletion in the experiments. We have chosen imparity of the wave functions but probably other  conditions compatible with $H^1$ functions and density depletion could be sought.

As $d$ is increased, the soliton loses stability. A further question would be to analyze the eigenfunctions of the linearized operator around the soliton solution, and in particular understand whether the first direction of instability is around a vortex/antivortex solution.

Related issues include the study of minimizers of the energy in the product spaces $\R^n \times M^k$, where $M^k$
is a compact Riemannian manifold under an $L^2$ constraint. It is proved in \cite{terra} that when the $L^2$ norm of the solution is sufficiently small, then the ground states
coincide with the corresponding $\R^n$ ground states, that is in our case, the soliton for a small torus. They also prove that above a critical mass the ground
states have nontrivial $M^k$ dependence.

The problem studied in this paper corresponds to what is known as zero velocity. Further works could be performed by minimizing the energy with fixed velocity or minimizing the energy with fixed momentum. The analysis of periodic solutions in the torus has been the topic of a very nice recent paper \cite{ruiz2}, where according to the size of the period, a mountain pass solution is constructed. Such methods could probably be generalized here in the case of the strip.

In the experimental papers \cite{taming} and \cite{statquantum}, in a one dimensional channel of polariton superfluids, the soliton breaks into arrays of vortex-antivortex. The mathematical analysis is probably similar though the equation there includes a forcing term and a dissipation term. The stationary one dimensional solutions have been studied in \cite{HPA} and the analysis of the the solitonic vortices in this context is an interesting open question.

 \section{Proof of Theorem \ref{theo1}}
 \subsection{Existence of a minimizer and asymptotic properties}
Using the soliton a a comparison function,  we find that the minimal energy is finite and is
bounded above by a constant times $d$.
From any  minimizing sequence, we may extract a subsequence which converges weakly in $H^1_{loc}$ to a limiting $u$ which satisfies the symmetry assumptions and therefore is a minimizer of (\ref{energy}) under conditions (\ref{eqsymmetry}). It also satisfies Neumann boundary conditions on the lines $\{y=\pm d\}$.

From the energy minimality, $|u|_\infty\le 1$ in the strip and by elliptic regularity, the gradient of $u$ is  uniformly bounded in the strip, and tends to $0$ as $x\to\pm\infty$. This first implies that $|u|\to 1$ as $x\to +\infty$: indeed, the gradient bound implies that if $|u(z)|<1-\eta$ for some $\eta>0$, then $|u|<1-\eta/2$ on a ball of radius $r_\eta>0$, and therefore the energy of $u$ on such a  ball is bounded below by some $c_\eta>0$. This cannot happen for arbitrarily large $z$ since the energy is finite. Hence $|u(x,y)|\to 1$ as $x\to +\infty$.

Moreover, since the gradient tends to $0$ in the infinite direction, we have that, as $x\to \pm\infty$,  $\|u(x,\cdot) - u(x,0)\|_{L^\infty(-d,d)}\to 0$. But conditions (\ref{eqsymmetry}) imply that $u(x,0)$ is purely imaginary. Since $|u(x,0)|$ tends to $1$ as $x\to\pm\infty$, it follows that for any $x$ such that $|x|$ is large enough depending on $\eta>0$, we have that $|u(x,y)-i|<\eta$ or $|u(x,y)+i| < \eta$. Therefore $u(x,y)$ converges to either $i$ or $-i$ as $x\to \pm\infty$ and, using (\ref{eqsymmetry}) again, that is $u(x,y) = \bar u(-x,y)$, the limits are opposite of one another.

\begin{rem}
 Let us point out that because of (\ref{eqsymmetry}) we have $u(0,0)=0$, therefore the energy minimum is not zero. Mathematically, the imparity condition (the second condition in (\ref{eqsymmetry})), can be viewed as a way to exclude an energy minimum of zero. It is an open question to choose a space of minimization including the first condition of (\ref{eqsymmetry}) with limits $i$ and $-i$ at plus and minus infinity, and such that $u$ vanishes somewhere.
\end{rem}

 \subsection{Uniqueness for small $d$} As a first step, we want to prove that when $d$ is small, the minimizer is purely imaginary, that is $u_1=\Re u\equiv 0$. We have, in the strip,
\begin{equation}\label{equ1}-\Delta u_1 + u_1 (|u |^2-1)=0,\end{equation}
together with Neumann boundary conditions on the lines $\{y = \pm d\}$. Because of (\ref{eqsym})-(\ref{eqsym2}), $u_1(x,0)$ has to be both odd and even so that we have $u_1(x,0)=0$. From  the Neumann boundary conditions, we have $\partial_y u_1(x,\pm d) = 0$.

We define $\tilde u_1(x,y) = u_1(x,y)$ for $y\in(0,d)$ and $\tilde u_1(x,y) = u_1(x,d-y)$ for $y\in[d,2d)$.   It follows that $\tilde u_1$ is  a solution  of (\ref{equ1}) on the strip $\R\times (0,2d)$ such that $\tilde u_1=0$ on the boundary $y=0$ and $y=2d$.

We know that $\tilde u_1$ tends to 0 at infinity since $u$ tends to $\pm i$. Moreover, it follows from the equation, as proved in \cite{bbh2}, that $|u|\leq 1$. Therefore, $\tilde u_1$ is a solution of $-\Delta \tilde u_1 + c(x) \tilde u_1=0$ with $c(x)$ uniformly bounded. For $d$ small enough, the Maximum Principle in narrow domains \cite{BN,BR,BNV} applies and implies that $\tilde u_1$ cannot reach a positive maximum or a negative minimum in $\R\times (0,2d)$: the proof consists in writing the equation satisfied by $w=\tilde u_1/g$ where $g= cos (\alpha y)$ where $\alpha=\pi /(8d)$ is chosen such that its square is bigger than the $L^\infty$ bound of $c(x)$. We have
$$-\Delta w-2 \frac {\nabla g} g \cdot \nabla w +(\alpha^2+c(x)) w=0.$$
  Since $\partial w/\partial n$ is zero on the boundary and $w$ tends to zero at infinity, and satisfies an equation for which the Maximum Principle holds, $w$ cannot reach a positive maximum or a negative minimum and is identically zero. Therefore $\tilde u_1$ is identically zero and  $u$ is purely imaginary.

Knowing that the solution is purely imaginary, its energy on each line $\{y\}\times \R$ can be bounded below by that of the soliton $i\tanh(x/\sqrt 2)$, which minimizes the energy among purely imaginary competitors satisfying the boundary conditions. Integrating with respect to $y$, we deduce that the energy of $u$ is bounded below by $2d$ times the energy of the soliton, and is bounded above by the same quantity since $u$ is a minimizer.
It follows that $\partial_y u$ is identically zero, and that the restriction of $u$ to any horizontal line is a (possibly translated) soliton. Since $\partial_y u = 0$, the translation is independent of $y$, and since $u(0,0)=0$, hence $u$ is the soliton.

\subsection{Upper bound for  the energy of a vortex for large $d$ and introduction of the solitonic vortex}
We now focus on the asymptotics of minimizers of $E_d$ for large $d$. We begin by computing an upper-bound for the minimal energy, for which a matching lower-bound will be proven below.

To prove the upper bound, we construct a test function which is an approximate  solution with a vortex at the origin. A first natural  test function is $u(z) = iz$ if $|z|\le 1$ and $u(z) = iz/|z|$ if $|z|>1$. It is not difficult to check that this yields the correct upper-bound but only up to a $O(1)$ term. Indeed, the integral on the strip, outside the unit disk is, in polar coordinates,
$$\int_{r = 1}^d \int_{-\frac \pi 2}^{\frac \pi 2}\frac1{r^2}\, r\,dr\, d\theta +\int_{r = d}^{+\infty}\int_{-\arcsin(d/r)}^{\arcsin(d/r)} \frac1{r^2}\, r\,dr\, d\theta.$$
The inner integral is equivalent to $2d/r^2$ as $r\to +\infty$, hence is convergent, but this does not provide the optimal constant.

In order to be more precise, we  construct a  test function with a vortex of degree one satisfying Neumann boundary conditions. This can be obtained by superimposing the known solution for a vortex without boundaries with those of image vortices. This provides a condition of no normal flow to the boundary for the gradient of the phase. While a single straight wall requires a single image vortex, the two parallel walls of the channel generate a doubly infinite array of image vortices.  The problem can be solved on the complex plane by introducing the meromorphic function $w$ having zeros (resp. poles) at the location of vortices of degree $+1$ (resp. $-1$) \cite{greengard,barton,toikka}. Then the phase is the argument of $w$ and the harmonic conjugate of the phase is the logarithm of its modulus. The advantage of considering the latter is that it has the same gradient energy as the phase but is single valued.

 The solution $w$ for the strip with a vortex at $z_0$ is given in \cite{greengard} for a strip $(0,H)$ and is
 $$\frac{\sinh\(\frac\pi{2H} (z - z_0)\)}{\sinh\(\frac\pi{2H} (z-\bar z_0)\)}.$$  It is found using a conformal mapping between the strip and the disc. Mapping the interval $(0,H)$ to $(-d,d)$ yields for the strip $\Omega_d$ with a vortex of degree $+1$ at $z = ia$
\begin{equation}\label{wa}
 w_a(z)=\frac{\sinh\(\frac\pi{4d} (z - ia)\)}{\sinh\(\frac\pi{4d} (z+ia+2id)\)} = i\frac{\sinh\(\frac\pi{4d} (z - ia)\)}{\cosh\(\frac\pi{4d} (z+ia)\)},
 \end{equation}
 so that, if $a = 0$,
 \begin{equation}\label{w0}
 w_0(z)=i\tanh\(\frac{\pi z}{4d}\).
 \end{equation}
 \begin{figure}[htbp]
\begin{center}\includegraphics[width=0.9\textwidth]{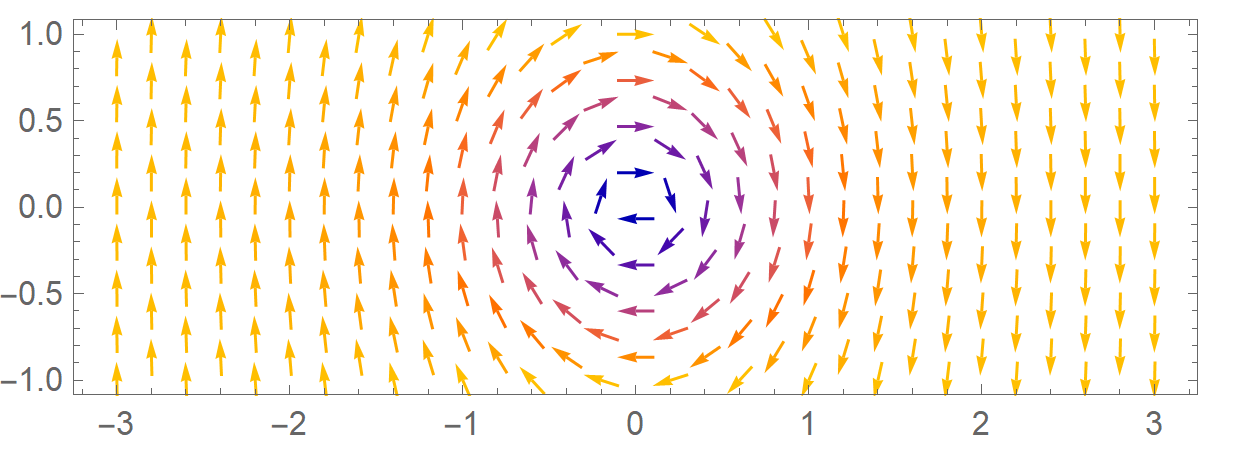}	
\caption{The function $w_0$ for a vortex at the origin.}\label{figw1}
\end{center}\end{figure}
 The function $w$ is plotted in Figure \ref{figw1}, and we see that it behaves as a soliton far away from the core.
 The phase $S$ of a function $-i \log \frac {\alpha + i\beta}{\gamma +i\delta}$ is such that $\tan S= \frac{\beta\gamma-\alpha\delta}{\alpha\gamma+\beta\gamma}$ so that the phase $\varphi_a$ is given from (\ref{wa}) by
 \begin{equation}\label{phia}\tan \varphi_a (x,y) =\frac  {\sinh (\frac {\pi x}{2d})\cos (\frac {\pi a}{2d} ) } {\sin (\frac {\pi y}{2d} )-\cosh (\frac {\pi x}{2d})\sin (\frac {\pi a}{2d} )}.\end{equation} This formula yields that the phase difference between $x=-\infty$ and $x=+\infty$ is related to the location of the vortex and is
 $$\pi-\pi \frac a d.$$ So in order to have a phase difference of $\pi$ (between $-i$ and $i$), the vortex has to be at the origin, that is $a=0$.
Then, the phase $\varphi_0$ and its harmonic conjugate $h_0$ are defined by
\begin{equation}\label{voho}\tan \varphi_0 (x,y) =\frac  {\sinh (\frac {\pi x}{2d}) } {\sin (\frac {\pi y}{4d} )}\quad
\hbox{ and }\end{equation} \begin{equation}\label{h0} h_0 (x,y)= \log\left|\tanh (\frac {\pi (x+iy)}{2d})\right|=\frac 12 \log \left |\frac{\cosh (\frac {\pi x}{2d})- \cos (\frac {\pi y}{2d})}{\cosh (\frac {\pi x}{2d})+\cos (\frac {\pi y}{2d})}
\right |.\end{equation}
The test function $u$ is now obtained as follows : first, we let $u(x,y) = v(x/d, y/d)$,  so that $v$ is defined in $\Omega_1$ instead of $\Omega_d$. We let $\ep = 1/d$ and choose $\rho >0$ which will be eventually chosen small, but large compared to $\ep$. We have, as noted earlier,
\begin{equation}\label{bbhenergy} E_d(u) = F_\ep(v) := \frac12\int_{\Omega_1}|\nabla v|^2 +\frac1{4\ep^2} \(1 - |v|^2\)^2.\end{equation}

In $D_\rho = \Omega_1\setminus B_\rho$, we define $v$ in agreement with the discussion above, that is $ v = e^{i\varphi_0}$, where $\varphi_0$  is defined in \eqref{voho}, with $d$ chosen equal to $1$. We thus have, in $D_\rho$,
$$|\nabla v|^2 = |\nabla \varphi_0|^2 = |\nabla h_0|^2.$$

To compute the energy of $v$ in $D_\rho$, we first note that the potential term in \eqref{bbhenergy} vanishes since $|v| = 1$ there, and then we integrate by parts as in \cite{bbh} to find that
$$\int_{D_\rho} |\nabla v|^2  = \int_{D_\rho} |\nabla h_0|^2 = \int_{\partial B_\rho} h_0\partial_\nu h_0 ,$$
where $\nu$ is the normal to the boundary pointing toward the center of the ball, and where we have used the fact, easily deduced from \eqref{h0}, that $h_0 = 0$ on $\partial \Omega_1$.
  Note also that, still from \eqref{w0}, the derivative of $h_0$  with respect to $x$ tends to $0$ at infinity exponentially fast so that there are no boundary terms at infinity.

From \eqref{w0} we have, for small $\rho$,
$$w_0(\rho e^{i\theta})  = i \frac{\pi}{4} \rho e^{i\theta}\(1+O(\rho^2)\).$$
Therefore, as $\rho\to 0$,
$$h_0(\rho e^{i\theta})\simeq\log (\frac {\pi\rho}{4}) ,\quad \partial_\nu h_0(\rho e^{i\theta})\simeq -\frac1\rho.$$
It follows after some calculations that, for $\rho$ small,
\begin{equation}\label{enrho}\frac 12 \int_{D_\rho} |\nabla h_0|^2 = - \pi\log (\frac {\pi\rho} {4}) +  o(1). \end{equation}
To define $v$ inside $B_\rho$, we note that, since $h_0$ is the harmonic conjugate   to $\varphi_0$ and since $\partial_\nu h_0 \simeq -1/\rho$ on $\partial B_\rho$, the phase $\phi_0$ is close to $\theta + c$ as $\rho\to 0$, for some constant $c$. We may then define $v$ in $B_\rho$ to be the minimizer of $F_\ep$ in $B_\rho$ with Dirichlet boundary condition $e^{i(\theta+c)}$ on the boundary, modified slightly to match the definition of $v$ outside $B_\rho$. Following \cite{bbh}, such a minimizer will have energy $\gamma + \pi\log(\rho/\ep)+o(1)$ if $\ep/\rho$ is small and thus, adding \eqref{enrho} we find.
\begin{equation}\label{upper} E_d(u) = F_\ep(v) = \pi\log\frac\rho\ep + \gamma - \pi\log\frac{\rho\pi} 4 +o(1) =  \pi\log d + \gamma - \pi\log\frac\pi 4 +o(1),\end{equation}
where $o(1)$ is small if both $\rho$ and $\ep/\rho$ are small. Both are true if $d$ is large and if we choose, for instance, $\rho = 1/\sqrt d$.

\subsection{Lower bound and proof of Theorem~\ref{theo1} completed}
Let $u_d$ be a minimizer for the energy on the strip $\Omega_d$. We have $E_d(u_d) \le E_d(u)$, where $u$ was constructed in the previous section and satisfies \eqref{upper}. We let as above $\ep = 1/d$ and $v_\ep(x,y) = u_d(dx,dy)$. Then $E_d(u_d) = F_\ep(v_\ep)$, so that
\begin{equation}\label{equpper}  F_\ep(v_\ep) \le  \pi\log\frac 1\ep + \gamma - \pi\log\frac\pi 4 +o(1).\end{equation}

We use the following covering property, which can be deduced straightforwardly from   (\cite{sandier}, \cite{jerrard}).

\begin{prop} Assume $v_\ep$ satisfies  (\ref{equpper}) and  that $|v|$ tends to $1$ as $x\to\pm\infty$. Then, for any $\ep < r < 1$,   the strip $\Omega_1 = \R\times (-1,1)$  may be covered by closed disjoints disks $D_1,\dots,D_k$  of total radius $r$,  such that $|v_\ep|>1/2$ outside the disks and such that the following lower bound holds on each disk:
\begin{equation}\label{eqlower} F_\ep(v_\ep,D_i)\ge \pi |n_i|\(\log \frac r\ep - C\).\end{equation}
Here $n_i(r)$ is the degree of the mapping $v_\ep/|v_\ep|$ from the circle $\partial D_i$ to the unit circle, or zero if $D_i$ is not contained in the strip.
\end{prop}

Note that, since $v_\ep(0,0) = u_d(0,0) =  0$, there is at least one disk and the origin belongs to one of the disks.

For any fixed $r>0$, comparing (\ref{equpper}) with  (\ref{eqlower}), we find that
$$\sum_i |n_i|\le \frac{\log(1/\ep) + C}{\log(r/\ep) - C}.$$
It follows that  $\sum_i |n_i|\le  1$  if $\ep$ is small enough (corresponding to $d$ large  enough), and then   that there is at most one disk with nonzero degree, and that, if it exists, its degree is either $+1$ or $-1$. We claim that
\begin{enumerate}
\item A disk, call it $D_0$,  such that the degree on $\partial D_0 = \pm 1$ exists.
\item $D_0$ must be the disk containing the origin.
\end{enumerate}

To prove the second claim we argue by contradiction. If $0\not\in D_0$, then $D_0$ is disjoint from $D_1 = \{-z\mid z\in D_0\}$. But the energy of $v_\ep$ on $D_1$ is the same as its energy on $D_0$, from (\ref{eqsymmetry}), therefore the total energy of $v_\ep$ is at least twice that on $D_0$, which contradicts the upper bound if $\ep$ is small enough.

To prove the first claim we argue again by contradiction. Assume all disks have degree zero, then choose a circle of radius $t$ centered at the origin and not intersecting the other disks --- this is possible if  $r<1$. Then  the degree $n$ of $v_\ep$ on $\{|z|=t\}$ must be zero. This is impossible because $n$ is odd: Indeed, denoting by $\phi$ the phase of $v_\ep$ and $\theta$ the polar angle on the circle, the conditions (\ref{eqsymmetry}) mean that, modulo $2\pi$, we have  $\phi(\pi-\theta) = -\phi(\theta)$ and $\phi (\theta+\pi) = \phi(\theta)+\pi.$ It follows that
$$\int_{-\pi}^0 \phi' = \int_0^\pi \phi',\quad \int_0^\pi \phi' = \pi +2k\pi,$$
and therefore the integral of $\phi'$ over the interval $(0,2\pi)$ belongs to $2\pi +4\pi\Z$ or, equivalently, the degree $n$ is odd. We assume below without loss of generality that the degree is $+1$

The rest of the lower bound follows the well-known line of arguments found in \cite{bbh}. Comparing again \eqref{equpper} and \eqref{eqlower}, we find that for any choice of $r$, the energy of $v_\ep$ outside a disk of radius smaller than $r$ which contains the origin remains bounded as $\ep\to 0$. Therefore, modulo a subsequence, $v_\ep$ converges outside the origin (weakly in $H^1_{\text{loc}}(\Omega_1\setminus\{0\}$ to be precise) to a map $v_0$ which must be of constant modulus $1$, must be harmonic, and must satisfy a Neumann boundary condition on $\partial \Omega_1$. These properties follow from the energy bound satisfied by $v_\ep$ and the fact that it is a minimizer, hence a critical point,  of $F_\ep$. The degree condition above further implies that in fact $v_0$ is the canonical harmonic map in $\Omega_1$ with Neumann boundary conditions and a single singularity of degree $+1$.  We deduce that $v_0 = e^{i\varphi_0}$ is exactly as above. The fact that the convergence holds  in fact in $C^1$ norm, locally in $\bar\Omega_1\setminus\{0\}$ is proved in \cite{bbh}, Theorem VI.1.

Then, for  any fixed $\rho>0$, we have
$$\liminf_{\ep\to 0} F_\ep(v_\ep,D_\rho)\ge \frac12 \int_{D_\rho} |\nabla \varphi_0|^2 =  - \pi\log (\frac {\pi\rho} {4}) +  o_\rho(1),$$
where $o_\rho(1)$ denotes a quantity which tends to $0$ when $\rho\to 0$.

On the other hand, since $v_\ep$ converges to $e^{i\varphi_0}$ on $\partial B_\rho$, its energy in $B_\rho$ is bounded below by that of the minimizer of $F_\ep$ on $B_\rho$ with boundary data $e^{i\varphi_0}$, which is equal to that of the minimizer with boundary condition $e^{i\theta}$, up to a $o_\rho(1)$ error. Thus
$$\liminf_{\ep\to 0} F_\ep(v_\ep,B_\rho)\ge \gamma + \pi\log\frac\rho\ep +  o_\rho(1).$$

Adding the two lower bounds we find
$$F_\ep(v_\ep)\ge  \pi\log\frac 1\ep + \gamma  -   \pi\log (\frac {\pi} {4}),$$
which matches \eqref{equpper} as $\ep\to 0$ and thus concludes the proof of Theorem~\ref{theo1}.

\section{Computation of the renormalized energy of a vortex-antivortex pair}
We have seen that for $d$ small, the ground state is a soliton, and for $d$ large, a solitonic vortex.
 In addition to the minimizers described above, and as we mentioned in introduction, experimental evidence and numerical simulations  indicate that in the case of intermediate $d$, the soliton solution destabilizes into a vortex-antivortex pair.  It is not clear how to prove the existence of such solutions, or wether they are the solution of some variational problem. However we are going to provide some computations  for a vortex-antivortex solution. A rigorous proof might follow the so-called ``gluing" approach initiated in \cite{pacard} for Ginzburg-Landau vortices or from the Moutain Pass solutions of \cite{ruiz1}.

 \begin{prop}When $d$ is large, there is a test function with a vortex at $(0,d/2)$ and an anti-vortex at $(0,-d/2)$ whose energy is asymptotically \begin{equation}\label{eqvav} 2\pi\log d- 2\pi\log (\frac {\pi} {2}) + 2\pi \log (\sin\frac{\pi a}d) +2\gamma+ o(1).\end{equation}\end{prop}

We compute the energy of an approximate solution $u_a$ in $\Omega_d=\R \times (-d,d) $, having a vortex located at $p_+ = (0,a)$ and an antivortex at $p_- = (0,-a)$. The assumption is that $d$ is large compared to the size of the vortex core, and that the distance between vortices as well as their distance to the boundary is also large with respect to the size of the core. Denoting by $\rho$ the core size, this means that
\begin{equation}\label{rho} \rho \ll a \hbox{ and } \rho \ll d-a.\end{equation}

 We approximate $u_a$ outside $B(p_+,\rho)\cup B(p_-,\rho)$ by a map of modulus one, i.e. we set $u_a = e^{i\varphi_a}$ there, with $\varphi_a$ denoting the phase, defined only modulo $2\pi$. Then the energy density for $u_a$ reduces to $|\nabla\varphi_a|^2$ and therefore, since we are looking for an energy minimizing solution with fixed vortices, we assume $\varphi_a$ is harmonic. Under the assumption \eqref{rho}, the energy inside the core may be assumed to be independent of the precise vortex location and is for two vortices
  \begin{equation}\label{eq2v}2\gamma +2 \pi\log(\rho)+o(1).\end{equation}
  Therefore we will now consider the energy outside the core as a function of $a$.

\begin{figure}[htbp]
\begin{center}\includegraphics[width=0.9\textwidth]{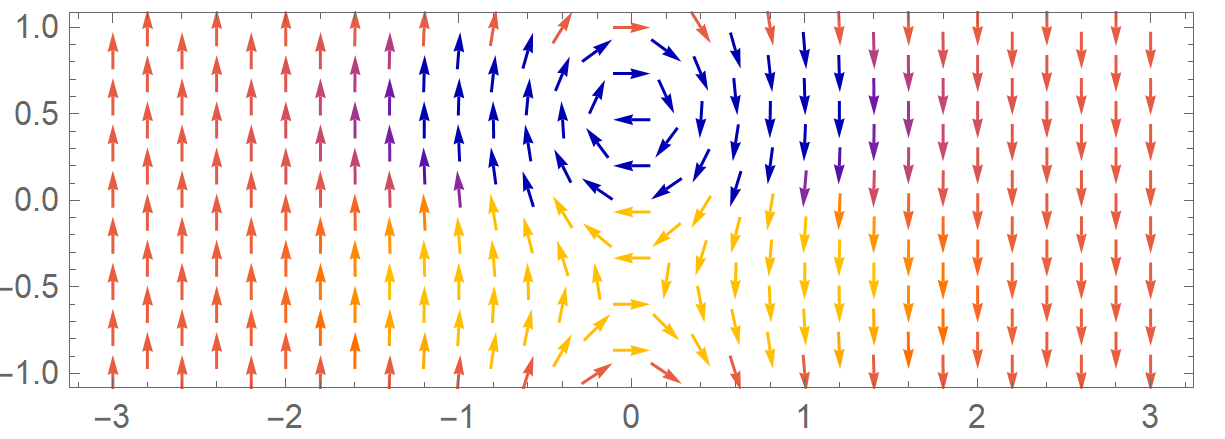}	
\caption{The function $W$ for a vortex/antivortex pair.}\label{figW2}
\end{center}\end{figure}
\begin{figure}[htbp]
\begin{center}\includegraphics[width=0.9\textwidth]{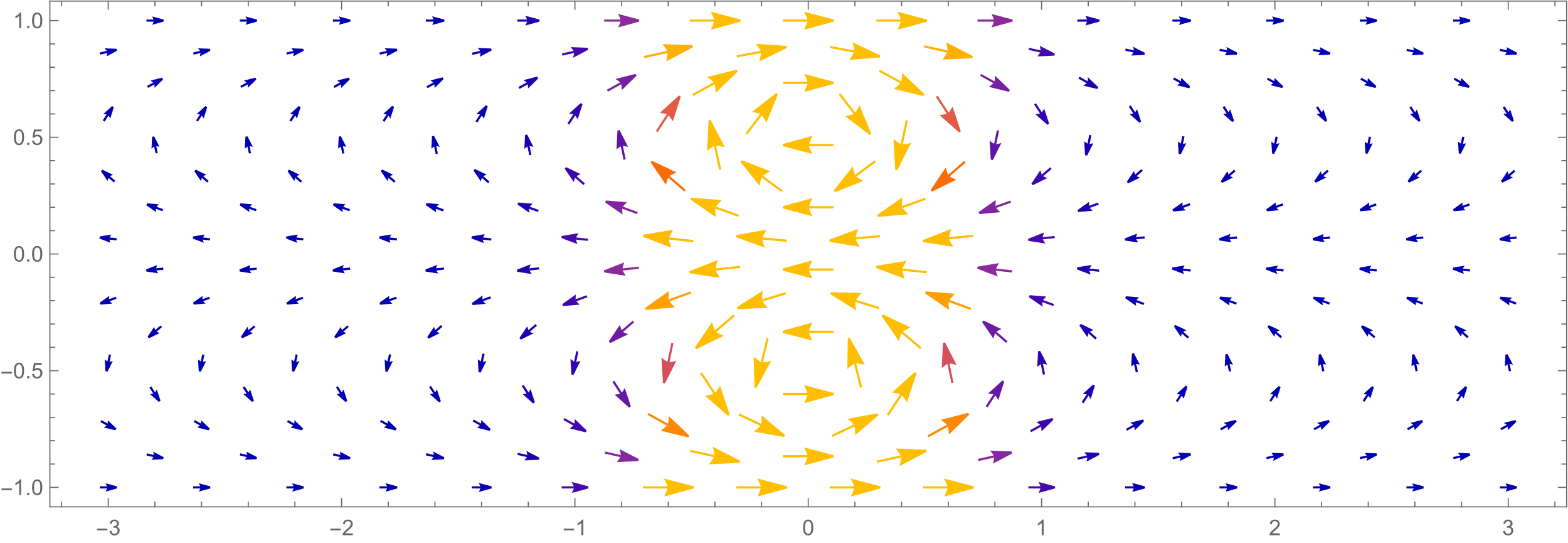}	
\caption{The gradient of the phase for a vortex, antivortex pair.}\label{figs1}
\end{center}\end{figure}
The harmonic conjugate  function to $\varphi_a$, call it $h_a$,  is single-valued. Since $\varphi_a$ satisfies a Neumann condition on $\partial \Omega_d$, the function $h_a$ satisfies a Dirichlet condition. In principle $h_a$ could take two different values on the top and bottom boundary, but this would give an infinite energy. Therefore the value is the same and $h_a$  solves the following equation:
 \begin{equation}\label{harmonic}
 \begin{cases} \Delta h_a = 2\pi (\delta_{p_+} - \delta_{p_-})  & \text{in $\Omega_d$}\\ h_a = 0 & \text{on $\partial\Omega_d$}\end{cases}\end{equation}
Note that the exact equation inside the vortex cores is unimportant for this calculation, since we wish to compute the energy away from the cores. This is why we introduce $\delta$-functions.

There is an explicit formula for $h_a$ that corresponds to single vortex in a double interval. We can in fact deduce it from \eqref{wa}: the meromorphic function corresponding to a vortex of degree $+1$ at $ia$ and a vortex of degree $-1$ at $-ia$ is $w_a/w_{-a}$, hence
$h_a(x,y) = \log|W(x+iy)|$, where
$$W(x+iy) = \frac{\sinh\(\frac\pi{2d} x\)\cos\(\frac\pi{2d} (y-a)\) + i \cosh\(\frac\pi{2d} x\)\sin\(\frac\pi{2d} (y-a)\)}{\sinh\(\frac\pi{2d} x\)\cos\(\frac\pi{2d} (y+a)\) + i \cosh\(\frac\pi{2d} x\)\sin\(\frac\pi{2d} (y+a)\)},$$
which may also be written as
\begin{equation}\label{w}W(z) = \frac{\sinh\(\frac\pi{2d} (z - ia)\)}{\sinh\(\frac\pi{2d} (z+ia)\)}.\end{equation} The function $W$ is plotted in Figure \ref{figW2}. Let us point out that the function $W$ is the same as the one for a single vortex in $a$ with $4d$ replaced by $2d$, which is consistent with the image vortices. The gradient of the phase is plotted in Figure \ref{figs1}.

To compute the energy outside the core we integrate by parts as in \cite{bbh} to find, writing
$$D_\rho = \Omega_d\setminus \(B(p_+,\rho)\cup B(p_-,\rho)\),$$
that
$$\int_{D_\rho} |\nabla h_a|^2 = \int_{\partial B(p_+,\rho)} h_a\partial_\nu h_a + \int_{\partial B(p_-,\rho)} h_a\partial_\nu h_a,$$
where $\nu$ is the normal to the boundary point toward the center of the ball. Note that  $h_a$ tends to 0 at infinity and its derivative with respect to $x$ tends to 0,  so that there are no terms at infinity.
 From \eqref{w} we have, for small $r$,
$$ W(i a + re^{i\theta})  = i \frac{\pi}{2d} re^{i\theta}\( \sin\frac{\pi a}d \)^{-1}+O(r^2),$$
 $$W(- i a + re^{i\theta})^{-1} = -  i \frac{\pi}{2d} re^{i\theta} \(\sin\frac{\pi a}d\)^{-1} +O(r^2).$$
Therefore, $$h_a(i a + re^{i\theta})=\log (\frac {\pi\rho} {2d}) - \log (\sin\frac{\pi a}d)$$ and on $\partial  B(p_+,\rho)$, $\partial_\nu h_a/\partial \nu= 1/\rho$.
It follows after some calculations that, for $\rho$ small,
\begin{equation}\label{hdrho}\frac 12 \int_{D_\rho} |\nabla h_a|^2 = - 2\pi\log (\frac {\pi\rho} {2d}) + 2\pi \log (\sin\frac{\pi a}d) + o(1). \end{equation}
 Let us point out that in the small $\rho$ approximation, this energy has a maximum for $a = d/2$, but on the other hand if we want to have the same structure as the soliton, that the limits are $\pm i$ at $\pm \infty$, the only value of $a$ is $a=d/2$. Indeed the formula for the phase is the same as (\ref{phia}) with $d$ replaced by $d/2$ and $a$ replaced by $d/2-a$. Note that in the whole space, the vortex anti-vortex solution does not have a phase difference between two infinite directions.

 We recollect (\ref{eq2v}) and (\ref{hdrho}) to find (\ref{eqvav}).

We see that when $d$ is large, (\ref{eqvav}) provides a higher energy than (\ref{bbh}), but it may be that for intermediate $d$, the vortex anti-vortex solution is the ground state.


\hfill

\noindent
{\bf Acknowledgements:} The authors would like to thank Frédéric Chevy for discussions on the physics of solitons and solitonic vortices. This project was supported by the grant CNRS-Prime TraDisQ1D.
\bibliography{ref}

\end{document}